\documentclass[12pt]{article}
\usepackage{amssymb}
\usepackage{latexsym,bm}
\usepackage{graphicx}
\usepackage{amsmath}
\usepackage{mathrsfs}
\usepackage{epstopdf}
\usepackage{tikz}
\usepackage{amsthm}
\usepackage{mathtools}
\usepackage{enumitem}

\date{}
\theoremstyle{remark}

\begin{document}
\title{The maximum number of paths of a given length in a nonhamiltonian graph}
\author{\hskip -10mm Chengli Li and  Xingzhi Zhan\thanks{Corresponding author.}}
\maketitle
\footnotetext[1]{Department of Mathematics,  Key Laboratory of MEA (Ministry of Education)
 \& Shanghai Key Laboratory of PMMP, East China Normal University, Shanghai 200241, China}
\footnotetext[2]{E-mail addresses: {lichengli0130@126.com (C. Li),}  {zhan@math.ecnu.edu.cn (X. Zhan).}}

\begin{abstract}
In 1980, Paul Erd\H{o}s posed the following problem: For every positive integer $n,$ determine a nonhamiltonian graph of order $n$ having the maximum number of Hamilton paths.
We solve the more general problem of determining the nonhamiltonian graphs of order $n$ having the maximum number of paths of length $k$ for given integers $n$ and $k$
with $1\le k\le n-1.$ The case $k=n-1$ gives a solution to Erd\H{o}s's problem and the case $k=1$ corresponds to a theorem due to Ore and Bondy.
\end{abstract}

{\bf Key words.} Hamilton path; nonhamiltonian graph; extremal graph; Bondy--Chv\'atal closure

{\bf Mathematics Subject Classification.} 05C30, 05C35, 05C38
\vskip 8mm

\section{Introduction}

We consider finite simple graphs and use standard terminology and notations from [4] and [8]. The {\it order} of a graph is its number of vertices, and the
{\it size} is its number of edges. For graphs we will use equality up to isomorphism, so $G=H$ means that $G$ and $H$ are isomorphic.
We denote by $K_n$ the complete graph of order $n.$ Following Bondy [2], we use the following notation throughout.

{\bf Notation 1.} Denote by $K_{n-1}\cdot K_2$ the graph of order $n$ obtained by identifying a vertex of $K_{n-1}$ with a vertex of $K_2.$

The graph $K_6\cdot K_2$ is depicted in Figure 1.
\begin{figure}[h]
\centering
\includegraphics[width=0.38\textwidth]{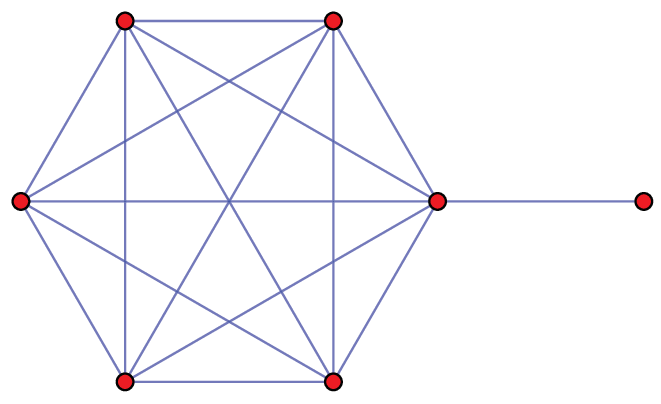}
\caption{The graph $K_6\cdot K_2$}
\end{figure}

The motivation for this research comes from a theorem due to Ore and Bondy, and a problem posed by Erd\H{o}s.

The following result on the maximum size of a nonhamiltonian graph of a given order is classic.

{\bf Theorem 1.} (Ore-Bondy) {\it The maximum size of a nonhamiltonian graph of order $n$ with $n\ge 6$ is $(n^2-3n+4)/2$ and this maximum size is attained uniquely
by the graph $K_{n-1}\cdot K_2.$ }

The extremal part of Theorem 1 was proved by Ore [7], and the structural part by Bondy [2]. For the case $n=5$ in Theorem 1, there are two extremal graphs.

In 1980, Paul Erd\H{o}s (see [5], p.201) posed the following problem.

{\bf Problem 1.} {\it For every positive integer $n,$ determine a nonhamiltonian graph of order $n$ having the maximum number of Hamilton paths.}

In this paper, we solve the more general problem of determining the nonhamiltonian graphs of order $n$ having the maximum number of paths of length $k$ for given integers $n$ and $k$ with $1\le k\le n-1.$ The case $k=n-1$ gives a solution to Problem 1 and the case $k=1$ corresponds to Theorem 1.

{\bf Notation 2.} For integers $n$ and $k$ with $n\ge k\ge 0,$ denote the permutation number by
$$
{\rm P}(n,k)=\begin{cases} n(n-1)(n-2)\cdots(n-k+1)\quad &{\rm if} \quad 1\le k\le n,\\
                           \qquad\qquad\quad\quad 1 \quad &{\rm if} \quad k=0.
             \end{cases}
$$

Our main result is as follows.

{\bf Theorem 2.} {\it Let $n$ and $k$ be integers with $n\ge 6 $ and $1\le k\le n-1.$ Then the maximum number of paths of length $k$ in a nonhamiltonian graph of order $n$ is
$$
{\rm P}(n-2,k-1)[(n-1)(n-1-k)/2+1]
$$
and this maximum number is attained uniquely by $K_{n-1}\cdot K_2.$}

The case $k=n-1$ of Theorem 2 gives a solution to Problem 1 as follows.

{\bf Corollary 3.} {\it Let $n$ be an integer with $n\ge 6.$ Then the maximum number of Hamilton paths in a nonhamiltonian graph of order $n$ is $(n-2)!$
and this maximum number is attained uniquely by $K_{n-1}\cdot K_2.$}

There are works in the literature counting the number of paths of a certain length in a graph ([1], [6], [9]).

\section{Proof}

We denote by $V(G)$ and $E(G)$ the vertex set and edge set of a graph $G,$ respectively, and denote by ${\rm deg}_G(x)$ the degree of a vertex $x$ in $G.$
Thus $e(G)\triangleq|E(G)|$ is the size of $G.$ We write $\overline{G}$ for the complement of $G$.

A nonhamiltonian graph $G$ is said to be {\it maximally nonhamiltonian} if for every non-edge $e\in E(\overline{G}),$ the graph $G+e$ is hamiltonian.
Let $G$ be a nonhamiltonian graph of order $n.$  Then every Hamilton cycle of $K_n$ contains at least
one edge of $\overline{G}$. For an integer $k$ with $1 \le k \le n-1,$ let $p_k(G)$ be the number of paths of length $k$ in $G$. In particular, $h(G)\triangleq p_{n-1}(G)$ is the number of Hamilton paths of $G$, and $p_1(G)=e(G).$

We record several easy facts in the following lemma.

{\bf Lemma 4.} {\it Consider the complete graph $K_n$ with $n \ge 4.$ For a given edge $ab,$ there exist exactly $(n-2)!$ Hamilton cycles containing $ab$ in $K_n.$
For two given adjacent edges $ab$ and $bc,$ there exist exactly $(n-3)!$ Hamilton cycles containing $ab$ and $bc$ in $K_n$.
For two given nonadjacent edges $ab$ and $cd,$ there exist exactly $2(n-3)!$ Hamilton cycles containing $ab$ and $cd$ in $K_n$.}

We shall use the Bondy--Chv\'atal closure theorem in the following form: if $u,v$
are nonadjacent vertices in a graph $G$ of order $n$ and $\deg_G(u)+\deg_G(v)\ge n$, then $G$ is
hamiltonian if and only if $G+uv$ is hamiltonian [3].

{\bf Lemma 5.} {\it Let $G$ be a maximally nonhamiltonian graph of order $n\ge 3.$  Then, for every
non-edge $uv\in E(\overline{G}),$ we have $\deg_{\overline{G}}(u)+\deg_{\overline{G}}(v)\ge n-1.$
}

{\bf Proof.}
Since $G$ is maximally nonhamiltonian, $G+uv$ is hamiltonian for every
$uv\in E(\overline{G})$.  If $\deg_G(u)+\deg_G(v)\ge n$, then the Bondy--Chv\'atal closure theorem would
imply that $G$ is hamiltonian, a contradiction.  Hence $\deg_G(u)+\deg_G(v)\le n-1$, and therefore
$\deg_{\overline{G}}(u)+\deg_{\overline{G}}(v)=2(n-1)-\bigl(\deg_G(u)+\deg_G(v)\bigr)\ge n-1$.
\hfill$\Box$

We first consider the case where $k=n-1$, and we require the following key lemma.

{\bf Lemma 6.}
{\it Let $G$ be a maximally nonhamiltonian graph of order $n \ge 4$ and let $m=e(\overline{G})$.  Then
$$h(G)\le \frac{m}{2m-n+2}(n-2)!.$$
In particular, if $m>n-2$, then $h(G)<(n-2)!$.
}

{\bf Proof.}
Let $\mathscr{C}$ be the set of Hamilton cycles of $K_n$.
For each integer $j \ge 1$, let
$$
\mathscr{C}_j=\{C\in\mathscr{C}: |E(C)\cap E(\overline{G})|=j\}~~\quad \text{and}~~
\quad
x_j=|\mathscr{C}_j|.
$$
Thus $x_j$ is the number of Hamilton cycles of $K_n$ containing exactly $j$ edges of $\overline{G}$.

Since $G$ is nonhamiltonian, every Hamilton cycle of $K_n$ uses at least one edge of
$\overline{G}$.  Therefore
\begin{equation}\label{eq:totalcycles}
\sum_{j\ge 1}x_j=|\mathscr{C}|=\frac{(n-1)!}{2}.
\end{equation}

We claim that $x_1=h(G).$ Indeed, if $C \in \mathscr{C}_1$, then deleting the unique edge of $E(C)\cap E(\overline{G})$ leaves a
Hamilton path of $G$. Conversely let $P$ be a Hamilton path of $G$ with endpoints $a$ and $b$.  Then $ab \notin E(G)$,
for otherwise $P+ab$ would be a Hamilton cycle of $G$, contradicting the assumption that $G$ is nonhamiltonian.  Hence $ab \in E(\overline{G})$, and $P+ab$ is a Hamilton cycle of $K_n$ containing exactly one edge of $\overline{G}$, namely $ab$.  This yields a bijection from the set $\mathscr{C}_1$ to the set of Hamilton paths of $G.$ Therefore, $x_1=h(G).$

Let
$$
\Omega=\{(C,e)|\, C\in \mathscr{C},\ e\in E(C)\cap E(\overline{G})\}.
$$
Counting $\Omega$ by first fixing the cycle $C$ gives $|\Omega|=\sum_{j\ge 1}j\,x_j$.
Counting $\Omega$ by first fixing the edge $e\in E(\overline{G})$ and then using
Lemma 4 gives $|\Omega|=m(n-2)!$.
Hence
\begin{equation}\label{eq:first-moment}
\sum_{j\ge 1}j\,x_j = m(n-2)!.
\end{equation}

Now define
$\beta\triangleq\sum_{j\ge 2}(j-1)x_j$.
Subtracting \eqref{eq:totalcycles} from \eqref{eq:first-moment}, we obtain
\begin{equation}\label{eq:Lvalue}
\beta=m(n-2)!-\frac{(n-1)!}{2}
=\frac{2m-n+1}{2}(n-2)!.
\end{equation}
By Theorem 1, we have $m\ge n-2$ and so
$2m-n+1\ge n-3>0$. It follows that $\beta>0.$

Let
$$
\mathscr{P}=
\{(C,\{e,f\})|\, C\in \mathscr{C},\ e,f\in E(C)\cap E(\overline{G}),\ e\neq f\},
$$
and let
$$
p\triangleq|\mathscr{P}|=\sum_{j\ge 2}\binom{j}{2}x_j.
$$

Let $s$ be the number of unordered pairs of adjacent edges of $\overline{G}$, and let $q$ be
 the number of unordered pairs of nonadjacent edges of $\overline{G}$.  So $s+q=\binom{m}{2}$.
Now count $\mathscr{P}$ by first fixing the pair of edges in $\overline{G}$.  By
Lemma 4, an adjacent pair belongs to exactly $(n-3)!$ Hamilton cycles of $K_n$, while a nonadjacent pair
belongs to exactly $2(n-3)!$ Hamilton cycles.  Hence
\begin{equation}\label{eq:Pcount}
p=s(n-3)!+2q(n-3)!=\bigl(m(m-1)-s\bigr)(n-3)!.
\end{equation}

We next obtain a lower bound on $s$.  Since each adjacent pair of edges of $\overline{G}$ has a
unique common endpoint,
$$
s=\sum_{v\in V(\overline{G})}\binom{\deg_{\overline{G}}(v)}{2}
=
\frac12\left(\sum_{v\in V(\overline{G})}\deg_{\overline{G}}(v)^2-2m\right).
$$
On the other hand,
$$
\sum_{uv\in E(\overline{G})}\bigl(\deg_{\overline{G}}(u)+\deg_{\overline{G}}(v)\bigr)
=
\sum_{v\in V(\overline{G})}\deg_{\overline{G}}(v)^2,
$$
because each vertex $v$ appears exactly once for each incident edge.  By
Lemma 5,
$$
\sum_{uv\in E(\overline{G})}\bigl(\deg_{\overline{G}}(u)+\deg_{\overline{G}}(v)\bigr)\ge m(n-1),
$$
and therefore
$\sum_{v\in V(\overline{G})}\deg_{\overline{G}}(v)^2\ge m(n-1)$.
Substituting into the expression for $s$ yields
$$
s\ge \frac12\bigl(m(n-1)-2m\bigr)=\frac{m(n-3)}2.
$$
Using this bound in \eqref{eq:Pcount}, we obtain
\begin{equation}\label{eq:Pupper}
p\le \frac{m(2m-n+1)}{2}(n-3)!.
\end{equation}

Let
$t\triangleq\sum_{j\ge 2}x_j$.
Then by \eqref{eq:totalcycles}, we have
$t=(n-1)!/2-x_1$.  We also have
$2p-\beta=\sum_{j\ge 2}(j-1)^2x_j.$
Applying Cauchy--Schwarz inequality to the sequences $\bigl(\sqrt{x_j}\bigr)_{j\ge 2}$ and
$\bigl((j-1)\sqrt{x_j}\bigr)_{j\ge 2}$ gives
$$
\left(\sum_{j\ge 2}(j-1)x_j\right)^2
\le
\left(\sum_{j\ge 2}x_j\right)\left(\sum_{j\ge 2}(j-1)^2x_j\right),
$$
that is,
$\beta^2 \le t(2p-\beta)$.
Since $\beta>0,$ we obtain $x_j>0$ for some $j\ge 2$
and so $2p-\beta>0.$
Therefore,
\begin{equation}\label{eq:t-lower}
t\ge \frac{\beta^2}{2p-\beta}.
\end{equation}

Combining \eqref{eq:Pupper} and \eqref{eq:Lvalue}, we find that
\begin{align}
2p-\beta
&\le
m(2m-n+1)(n-3)!
-
\frac{2m-n+1}{2}(n-2)! \notag\\
&=
\frac{(2m-n+1)(2m-n+2)}{2}(n-3)!.
\label{eq:2PminusL-upper}
\end{align}
Substituting \eqref{eq:Lvalue} and \eqref{eq:2PminusL-upper} into \eqref{eq:t-lower} gives
$$
t
\ge
\frac{
\left(\frac{2m-n+1}{2}(n-2)!\right)^2
}{
\frac{(2m-n+1)(2m-n+2)}{2}(n-3)!
}
=
\frac{(2m-n+1)(n-2)}{2(2m-n+2)}(n-2)!.
$$
Finally, using $h(G)=x_1$ and \eqref{eq:totalcycles}, we obtain
\begin{align*}
h(G)=x_1
&=
\frac{(n-1)!}{2}-t\\
&\le
\frac{n-1}{2}(n-2)!
-
\frac{(2m-n+1)(n-2)}{2(2m-n+2)}(n-2)!\\
&=
\frac{m}{2m-n+2}(n-2)!.
\end{align*}
This proves the first assertion.  If $m>n-2$, then
$\frac{m}{2m-n+2}<1$,
so $h(G)<(n-2)!$.
\hfill$\Box$

Now we finish the Hamilton path case.

{\bf Lemma 7.} {\it The maximum number of Hamilton paths in a nonhamiltonian graph of order $n\ge 6$ is $(n-2)!$, and this maximum number is attained uniquely by $K_{n-1}\cdot K_2.$
}

{\bf Proof.} Let $G$ be any nonhamiltonian graph of order $n$ attaining the maximum number of Hamilton paths. Adding edges if necessary, it is clear that there exists a maximally
nonhamiltonian graph $H$ such that $G$ is a spanning subgraph of $H.$ Continue using the natation $h(R)$ to denote the number of Hamilton paths in a graph $R.$ Since
$h(K_{n-1}\cdot K_2)=(n-2)!,$ we have $h(H)\ge h(G)\ge (n-2)!.$ By Lemma 6, we deduce that $e(\overline{H})\le n-2;$ i.e., $e(H)\ge (n^2-3n+4)/2.$ By Theorem 1, $H=K_{n-1}\cdot K_2.$
Since every edge of $K_{n-1}\cdot K_2$ lies in at least one Hamilton path, $G=H=K_{n-1}\cdot K_2.$ \hfill$\Box$

Next we prove the main theorem.

{\bf Proof of Theorem 2.}
The case $k=n-1$ is exactly Lemma 7.  It therefore remains to consider $1\le k \le n-2$.

Let $G$ be a nonhamiltonian graph of order $n.$  As in the proof of Lemma 6, let $\mathscr{C}$ be the set of Hamilton cycles of $K_n$, and for $j\ge 1$, let
$x_j=|\{C\in \mathscr{C}: |E(C)\cap E(\overline{G})|=j\}|$.
Then
\begin{equation}\label{eq:section4-total}
\sum_{j\ge 1}x_j=\frac{(n-1)!}{2}
\qquad\text{and}\qquad
x_1=h(G).
\end{equation}

For a Hamilton cycle $C\in \mathscr{C}$, let $s_k(C)$ be the number of paths of length $k$ in both $G$ and $C.$  We first show that
\begin{equation}\label{eq:segment-double-count}
(n-k-1)!\,p_k(G)=\sum_{C\in \mathscr{C}}s_k(C).
\end{equation}

Fix a path
$P=v_0v_1\cdots v_k$ of length $k$ in $G$.   A Hamilton cycle of $K_n$ contains $P$  if and only if it has a  form
$v_0,v_1,\dots,v_k,w_1,\dots,w_{n-k-1},v_0,$
where $(w_1,\dots,w_{n-k-1})$ is a permutation of the remaining vertices.
Thus exactly $(n-k-1)!$ Hamilton cycles of $K_n$ contain $P.$ Summing over all paths of length $k$ in $G$ proves \eqref{eq:segment-double-count}.

Next we bound $s_k(C)$ in terms of the number of edges of $\overline{G}$ lying in $C$.
Suppose that $C$ contains exactly one edge of $\overline{G}$.
Deleting that edge leaves a Hamilton path of $G$ on $n$ vertices, and such a path has exactly
$n-k$ consecutive subpaths of length $k$.  Hence
$s_k(C)=n-k$.

Suppose that $C$ contains $j\ge 2$ edges of $\overline{G}$.
Then the edges of $C$ that belong to $G$ split into $j$ paths along the cycle.
Let their lengths be $a_1,\dots,a_j$, so each $a_i\ge 0$ and
$a_1+\cdots+a_j=n-j$.
A path of length $a_i$ contains exactly $\max\{a_i-k+1,0\}$ consecutive subpaths of length
$k$.  Therefore
$$s_k(C)=\sum_{i=1}^j \max\{a_i-k+1,0\}.$$
Let
$I=\{i\in\{1,\dots,j\}|\,  a_i\ge k\}$.
If $I=\emptyset$, then $s_k(C)=0$.  If
$I\neq\emptyset$, then by the fact that $j\ge 2$, we have
\begin{align*}
s_k(C)=\sum_{i\in I}(a_i-k+1)
&=
\sum_{i\in I}a_i-|I|(k-1)\\
&\le
\sum_{i=1}^j a_i-(k-1)\\
&=
 n-j-k+1\\
&\le
 n-k-1.
\end{align*}

Using \eqref{eq:segment-double-count}, we conclude that
\begin{align}
(n-k-1)!\,p_k(G)=
\sum_{C\in \mathscr{C}}s_k(C)
&\le
(n-k)x_1+(n-k-1)\sum_{j\ge 2}x_j \notag\\
&=
(n-k-1)\sum_{j\ge 1}x_j+x_1 \notag\\
&=
(n-k-1)\frac{(n-1)!}{2}+x_1.
\label{eq:main-ineq}
\end{align}
By Lemma 7,
$x_1=h(G)\le (n-2)!$.
Substituting this into \eqref{eq:main-ineq} yields
$$
(n-k-1)!\,p_k(G)
\le
(n-k-1)\frac{(n-1)!}{2}+(n-2)!.
$$
Dividing by $(n-k-1)!$ gives
\begin{align}
p_k(G)&\le \frac{(n-1)!}{2(n-k-2)!}+\frac{(n-2)!}{(n-k-1)!}\notag\\
      &={\rm P}(n-2,k-1)[(n-1)(n-k-1)/2+1],
\end{align}
which proves the upper bound.

One readily verifies that $K_{n-1}\cdot K_2$ attains the upper bound in (11).
It remains to prove the uniqueness.  Let $G$ be a nonhamiltonian graph of order $n\ge 6$
with $p_k(G)=p_k(K_{n-1}\cdot K_2)$.
Since
$$
(n-k-1)!\,p_k(K_{n-1}\cdot K_2)
=
(n-k-1)\frac{(n-1)!}{2}+(n-2)!,
$$
the equality $p_k(G)=p_k(K_{n-1}\cdot K_2)$ together with \eqref{eq:main-ineq} implies that
$$
(n-k-1)\frac{(n-1)!}{2}+(n-2)!
\le
(n-k-1)\frac{(n-1)!}{2}+x_1.
$$
Hence
$h(G)=x_1\ge (n-2)!$. By Lemma 7, we have $G= K_{n-1}\cdot K_2$. This completes the proof of Theorem 2. \hfill$\Box$

\vskip 5mm
{\bf Acknowledgement.} This research  was supported by the NSFC grant 12271170 and Science and Technology Commission of Shanghai Municipality grant 22DZ2229014.

\end{document}